\documentclass[12pt]{article}
\newtheorem{theorem}{Theorem}[section]
\newtheorem{lemma}[theorem]{Lemma}
\newtheorem{cor}[theorem]{Corollary}
\newtheorem{prop}[theorem]{Proposition}
\newtheorem{conj}[theorem]{Conjecture}

\usepackage{latexsym,amssymb}
\usepackage{amsmath}

\newcommand{\Z}{\hbox{\bf Z}}

\newcommand{\qfd}{\hfill $\Box$}
\newcommand{\F}{\mathcal{F}}
\newcommand{\Sum}[2]{\underset{#1}{\overset{#2}{\sum}}}
\newcommand{\be}{\begin{equation}}
\newcommand{\ee}{\end{equation}}
\newcommand{\ol}[1]{\overline{#1}}

\begin{document}
\title{{Weighted Sequences in Finite Cyclic Groups}}
\author{{David J. Grynkiewicz}\thanks{Departament de Matematica Aplicada i Telematica,
  Universitat Politecnica de Catalunya,
  Campus Nord Edifici C3,
  Jordi Girona Salgado 1-3,
  08034-Barcelona, Spain. email: diambri@hotmail.com.
  Supported in part by the National Science Foundation, as an
MPS-DRF postdoctoral fellow, under grant DMS-0502193} and {Jujuan
Zhuang }\thanks{Institute of Mathematics, Dalian Maritime
University, Dalian, 116024, P.R. China. Corresponding email:
jjzhuang1979@yahoo.com.cn.}}

\maketitle
\begin{abstract}
Let $p>7$ be a prime, let $G=\Z/p\Z$, and let $S_1=\prod_{i=1}^p
g_i$ and $S_2=\prod_{i=1}^p h_i$ be two sequences with terms from
$G$. Suppose that the maximum multiplicity of a term from either
$S_1$ or $S_2$ is at most $\frac{2p+1}{5}$. Then we show that, for
each $g\in G$, there exists a permutation $\sigma$ of $1,2,\ldots,
p$ such that $g=\Sum{i=1}{p}(g_i\cdot h_{\sigma(i)})$. The question
is related to a conjecture of A. Bialostocki concerning weighted
subsequence sums and the Erd\H{o}s-Ginzburg-Ziv Theorem.
\end{abstract}

\section{Introduction}

\indent \indent Let $G$ be a finite abelian group (written
additively), and let $\F(G)$ denote the free abelian monoid over $G$
with monoid operation written multiplicatively and given by
concatenation, i.e., $\F(G)$ consists of all multi-sets over $G$,
and an element $S\in\F(G)$, which we refer to as a \emph{sequence},
is written in the form $S=\prod_{i=1}^kg_i=\prod_{g\in G}g^{v_g(S)}$
with $g_i\in G$, where $v_g(S)\in \mathbb{N}_0$ is the {\sl
multiplicity of $g$ in $S$} and $k$ is {\sl the length of} $S$,
denoted by $|S|=k$. Set
$$h(S)=\max_{g\in G}\{v_g(S)\}.$$ If $h(S)\leq
1$, then we call $S$ a square-free sequence in $G$, in which case we
may also regard $S$ as a subset of $G$. A sequence $T$ is {\sl a
subsequence of $S$}, which we denote by $T|S$, if $v_g(T)\leq
v_g(S)$ for every $g\in G$. By $\sigma(S)$ we denote the sum of all
terms in $S=\prod_{i=1}^kg_i$, that is $\sigma(S)= \Sum{i=1}{k}g_i$.
If $G$ is also a ring with multiplicative operation denoted $\cdot$, and  $S,\,T\in \F(G)$ with
$S=\prod_{i=1}^{k_1}g_i$, $T=\prod_{i=1}^{k_2}h_i$ and
$r=\min\{k_1,\,k_2\}$, then define
$$S\cdot T:=\{\Sum{i=1}{r}(g_{\sigma_1(i)}\cdot h_{\sigma_2(i)})\mid
\sigma_i\mbox{ a permutation of } 1,2,\ldots,k_i,\,i=1,2\}.$$

In 1996, Y. Caro made the following conjecture \cite{Caro}, which
can be regarded as a weighted version of the Erd\H{o}s-Ginzburg-Ziv
Theorem \cite{egz} (which is the case $W=1^n$).

\begin{conj}\label{conj:1.1}
Let $G=\Z/n\Z$. If $S,\,W\in \F(G)$ with $|W|=k$, $\sigma(W)=0$ and
$|S|=n+k-1$, then $0\in W\cdot S$.
\end{conj}

For prime cyclic groups, Conjecture \ref{conj:1.1} was confirmed by
N. Alon, A. Bialostocki and Y. Caro (see \cite{Caro}). W. Gao and X.
Jin showed, in particular, that Conjecture \ref{conj:1.1} is true if
$n=p^2$ for some prime $p\in \mathbb{P}$ (see \cite{GJ}), and more
recently, a complete confirmation of Conjecture \ref{conj:1.1} was
found by D. Grynkiewicz (see \cite{Gryn1}).

On the basis of  Conjecture \ref{conj:1.1}, A. Bialostocki made the
following conjecture \cite{bialostocki}.

\begin{conj}\label{conj:1.2}
Let $G=\Z/n\Z$ with $n$ even. If $S_1,\,S_2\in \F(G)$ with
$|S_1|=|S_2|=n$ and $\sigma(S_1)=\sigma(S_2)=0$, then $0\in S_1\cdot
S_2$.
\end{conj}

The example \be\label{bial-example}S_1=
0^{n-2}1(-1),~S_2=012\cdots(n-1)\ee was given to show Conjecture
\ref{conj:1.2} could not hold for odd $n$. Additionally, A.
Bialostocki confirmed Conjecture \ref{conj:1.2} for small numbers
using a computer.

A related question is to ask what conditions guarantee that
$S_1\cdot S_2=G$, so that every element of $G$, including zero, can
be represented as a sum of products between the terms of $S_1$ and
$S_2$. Of course, if $S_1$ has only one distinct term, then
$|S_1\cdot S_2|=1$, so some condition, say either on the
multiplicity of terms or on the number of distinct terms, is indeed
needed.

In this paper, we show, for a group of prime order $p>3$,  that
$h(S_i)$ being small is enough to guarantee $S_1\cdot S_2=G$; note
that $S_1\cdot S_2=G$ implies $0\in S_1\cdot S_2$ as a particular
consequence. Our main result is the following.

\begin{theorem}\label{theorem:1.3}
Let $p>3$ be a prime, let $G=\Z/p\Z$, and let $S_1,\,S_2\in\F(G)$
with $|S_1|=|S_2|=p$. If $p\neq 7$ and $\max\{h(S_1), h(S_2)\}\leq
\frac{2p+1}{5}$, or $p=7$ and $\max\{h(S_1), h(S_2)\}\leq 2$, then
$S_1\cdot S_2=G$.
\end{theorem}

Let $G=\Z/n\Z$. If $n\equiv -1\mod 4$, then the example
\be\label{new-bial-example}S_1=S_2=
0^{\frac{n-1}{2}}1^{\frac{n-1}{2}}(\frac{n+1}{2})\ee has
$\sigma(S_1)=\sigma(S_2)=0$, $\max\{h(S_1),\,h(S_2)\}=\frac{n-1}{2}$
and $0\notin S_1\cdot S_2$, giving an additional counterexample to
the possibility of Conjecture \ref{conj:1.2} holding for odd order
groups. It also shows that the bound $\max\{h(S_1), h(S_2)\}\leq 2$
for $p=7$ is tight in Theorem \ref{theorem:1.3}. The example
$S_1=S_2=0^31^2$ shows that the bound for $p=5$ is tight, and the
example given in (\ref{bial-example}) for $n=3$ shows that the
theorem cannot hold for $p=3$. Finally, letting $x=
\lceil\frac{2n+2}{5}\rceil$, the example
\be\label{newer-bial-example}S_1=S_2=0^x1^x2^{n-2x},\ee for $n>7$,
has $\max\{h(S_1),\,h(S_2)\}=\lceil\frac{2n+2}{5}\rceil$ and
$S_1\cdot S_2\subseteq [n-2x,4(n-2x)+x]$, so that $|S_1\cdot
S_2|\leq 3n-5x+1\leq n-1$, showing that the bound $\frac{2p+1}{5}$
from Theorem \ref{theorem:1.3} is also best possible.

\section{The Proof of the Main Result}

\indent \indent To prove Theorem \ref{theorem:1.3}, we need some
preliminaries. Given subsets $A$ and $B$ of an abelian group $G$,
their sumset is the set of all possible pairwise sums:
$$A+B:=\{a+b\mid a\in A,\,b\in B\}.$$ We use $\ol{A}$ to denote the complement of $A$ in $G$.
For a prime order group, we have the following classical inequality
\cite{cdt}.

\begin{theorem}[Cauchy-Davenport Theorem]\label{lemma:1.4} Let $p$ be a prime,
if $A$ and $B$ are nonempty subsets of $\Z/p\Z,$ then
$$
|A+B| \geq \min\{p, |A|+|B|-1\}.$$
\end{theorem}

The case when equality holds in the Cauchy-Davenport bound was
addressed by A. Vosper \cite{Vosper}.

\begin{theorem}[Vosper's Theorem]\label{vosper-lemma:1.4} Let $p$ be a prime,
and let $A,\,B\subseteq \Z/p\Z$ with $|A|,\,|B|\geq 2$. If
$$|A+B|=|A|+|B|-1\leq p-2,$$ then $A$ and $B$ are arithmetic
progressions with common difference. If $$|A+B|=|A|+|B|-1=p-1,$$
then $A=x-\ol{B}$ for some $x\in \Z/p\Z$.
\end{theorem}

As an immediate corollary of Theorem \ref{vosper-lemma:1.4}, we have
the following.

\begin{cor}\label{vosper-cor} Let $p$ be a prime,
and let $A,\,B\subseteq \Z/p\Z$ with $|A|,\,|B|\geq 2$. If
$$|A+B|=|A|+|B|-1<p$$ and $B$ is an arithmetic progression with difference $d$,
then $A$ and $A+B$ are also arithmetic progressions with difference
$d$.
\end{cor}

We will also need the following basic proposition (an immediate
consequence of Lemma 1 in \cite{ham}).

\begin{prop}\label{propp} Let $p$ be a prime,
and let $A,\,B\subseteq \Z/p\Z$ be nonempty with $|B|\geq 3$. If
$$|A+B|\leq |A|+|B|< p,$$ and $B$ is an arithmetic progression with difference $d$,
then $A+B$ is also an arithmetic progression with difference $d$.
\end{prop}

The following will be the key lemma used in the proof of Theorem
\ref{theorem:1.3}.

\begin{lemma}\label{lemma:1.5}  Let $p>3$ be a prime, let $G=\Z/p\Z$, and let
$U,\,V\in \F(G)$ be square-free with $|U|=|V|=3$. Then $|U\cdot
V|\geq 4;$ furthermore, assuming $p>7$, then equality is only
possible if $U\cdot V$ is not an arithmetic progression and either
$U$ and $V$ are both arithmetic progressions or else, up to affine
transformation, $U=01x$ and $V=01y$ with $x$ and $y$ the two
distinct roots of $z^2-z+1$.
\end{lemma}

\begin{proof}
By an appropriate pair of affine transformations (maps of the form
$z\mapsto \alpha z+\beta$ with $\alpha,\,\beta\in \Z/p\Z$ and
$\alpha\neq 0$), we may w.l.o.g. assume $U=01x$ and $V=01y$ with
$x,\,y\notin\{0,1\}$. By possibly applying the affine transformation
$z\mapsto -z+1$ to $U$, we may assume $x\neq y$ unless
$x=y=\frac{p+1}{2}$, in which case we may instead assume $x\neq y$
by applying the affine transformation $z\mapsto 2z$ to $U$ and the
affine transformation $z\mapsto -2z+1$ to $V$. Observe that
\be\label{UV-sumset}U\cdot V=\{1,x,y,xy,x+y,xy+1\}.\ee Hence, since
$x,\,y\notin\{0,1\}$ and $x\neq y$, it follows that $\{1,x,y\}$ is a
cardinality $3$ subset of $U\cdot V$. Consequently, if $|U\cdot
V|<4$, then $xy,\,x+y,\,xy+1\in \{1,x,y\}$.

Since $x,\,y\notin\{0,1\}$, if $xy\in \{1,x,y\}$, then
\be\label{one}xy=1;\ee if $x+y\in \{1,x,y\}$, then \be\label{two}
x+y=1;\ee and if $xy+1\in \{1,x,y\}$, then w.l.o.g. \be\label{three}
xy=x-1.\ee

If (\ref{one}), (\ref{two}) and (\ref{three}) all hold, then
(\ref{one}) and (\ref{three}) imply $x=2$, whence (\ref{two})
implies $y=-1$; thus (\ref{one}) implies $3=0$, contradicting that
$p>3$. As a result, we conclude that $|U\cdot V|\geq 4$. To prove
the second part of the lemma, we now assume $p>7$ and $|U\cdot
V|=4$.

First suppose that either $U$ or $V$ is an arithmetic progression,
say $U$. Thus w.l.o.g. $U=01(-1)$. If $V$ is also an arithmetic
progression, then by an appropriate affine transformation we obtain
$V=012$; consequently, $U\cdot V=\{\pm 1,\pm 2\}$, whence $U\cdot
V+U\cdot V=\pm\{0,1,2,3,4\}$, so that $|U\cdot V+U\cdot V|=9>|U\cdot
V|+|U\cdot V|-1$ for $p>7$, which implies $U\cdot V$ is not an
arithmetic progression, as desired. Therefore we may instead assume
$V$ is not an arithmetic progression.

Since $U=01(-1)$, $y\neq x$ and $y\notin\{0,1\}$, it follows in view
of (\ref{UV-sumset}) that $\{\pm 1,\pm y\}$ is a cardinality $4$
subset of $U\cdot V$. Hence, since $|U\cdot V|=4$, it follows that
$$U\cdot V=\{\pm 1,\pm y\}.$$ Thus it
follows in view of (\ref{UV-sumset}) that $y-1\in \{\pm 1,\pm y\}$,
whence $y\notin \{0,1\}$ implies either $y=2$ or $y=\frac{p+1}{2}$.
However, in either case $V$ is an arithmetic progression, contrary
to assumption. So it remains to handle the case when neither $U$ nor
$V$ is an arithmetic progression, and hence $x,\,y\notin
\{-1,0,1,2,\frac{p+1}{2}\}$.

If (\ref{one}) and (\ref{three}) hold, then $x=2$, while if
(\ref{two}) and (\ref{three}) hold, then $y\notin \{0,1\}$ implies
$x=-1$. Both cases contradict that $x\notin
\{-1,0,1,2,\frac{p+1}{2}\}$.

Suppose (\ref{one}) and (\ref{two}) hold. Then $x\neq y$ implies
that $x$ and $y$ are the two distinct roots of $z^2-z+1$. Moreover,
(\ref{UV-sumset}) implies that $U\cdot V=\{1,2,x,y\}$, whence
(\ref{two}) gives \be\label{UV-take3} \{1,2,3,4,
x+1,x+2,y+1,y+2\}\subseteq U\cdot V+U\cdot V.\ee Suppose $|U\cdot
V+U\cdot V|\leq 7$. Then at least one of the following cases holds:
$x\in \{-1,0,1,2,3\}$, $y\in \{-1,0,1,2,3\}$, $x=y$, $x=y+1$ or
$y=x+1$. If $x=y+1$ or $y=x+1$, say $x=y+1$, then (\ref{two})
implies $y=0$, a contradiction. Consequently, since $x,\,y\notin
\{-1,0,1,2,\frac{p+1}{2}\}$ and $x\neq y$, we conclude that either
$x=3$ or $y=3$, say $y=3$.  Hence (\ref{two}) implies $x=-2$, whence
(\ref{one}) yields $7=0$, contradicting that $p>7$. So we conclude
that $|U\cdot V+U\cdot V|\geq 8>|U\cdot V|+|U\cdot V|-1$, and thus
that $U\cdot V$ cannot be an arithmetic progression, as desired.

From the previous two paragraphs, we conclude that at most one of
(\ref{one}), (\ref{two}) and (\ref{three}) can hold, and thus that
at least two of the quantities $xy$, $x+y$ and $xy+1$ are not
contained in $\{1,x,y\}$. Since $|U\cdot V|=4$, all the quantities
$xy$, $x+y$ and $xy+1$ not contained in $\{1,x,y\}$ must be equal.
Thus, as $xy+1\neq xy$, we see that $\{1,x,y\}\cap \{xy,x+y,xy+1\}$
is nonempty. Hence, if $xy=x+y$ are the two quantities outside
$\{1,x,y\}$, then (\ref{three}) holds and so $x+y=xy=x-1$,
contradicting that $y\notin \{-1,0,1,2,\frac{p+1}{2}\}$; while if
$x+y=xy+1$ are the two quantities outside $\{1,x,y\}$, then
(\ref{one}) holds and so $x+y=xy+1=2$, which when combined with
(\ref{one}) yields $(y-1)^2=0$, contradicting that $y\notin
\{0,1\}$. This completes the proof.\qfd
\end{proof}

We now proceed with the proof of Theorem \ref{theorem:1.3}.

\begin{proof} First suppose $p>7$ (we will handle the cases $p\leq 7$ afterwards), and
assume by contradiction that $|S_1\cdot S_2|\leq p-1$. Let
$h:=\lfloor\frac{2p+1}{5}\rfloor$ and let $S_1=U_1\cdots U_h$ and
$S_2=V_1\cdots V_h$ be factorizations of $S_1$ and $S_2$ into
square-free subsequences $U_i$ and $V_i$ such that $|U_i|=|V_i|=3$
for $i\leq p-2h$ and $|U_i|=|V_i|=2$ for $i>p-2h$ (such
factorizations are easily seen to exist in view of $\max\{h(A),
h(B)\}\leq h$, $\frac{p}{3}<h<\frac{p}{2}$ and $|S_1|=|S_2|=p$; see
for instance \cite{acta}). Let $A_i:=U_i\cdot V_i$ for
$i=1,\ldots,h$. Note $|A_i|=2$ for $i\geq p-2h+1$ and that Lemma
\ref{lemma:1.5} implies $|A_i|\geq 4$ for $i\leq p-2h$, with
$|A_i|=4$ possible only if $A_i$ is not an arithmetic progression.
Also, $\Sum{i=1}{h}A_i\subseteq S_1\cdot S_2$, and hence $|S_1\cdot
S_2|\leq p-1$ implies \be\label{notallp}|\Sum{i=1}{h}A_i|\leq
p-1.\ee Thus, in view of Theorem \ref{lemma:1.4} applied to
$\Sum{i=2}{h}A_i+A_1$, we conclude that
\be\label{notallp-partii}|\Sum{i=2}{h}A_i|\leq p-4.\ee

Since $h>\frac{p}{3}$ (in view of $p\geq 11$), it follows that
$|A_h|=2$, and thus $A_h$ is an arithmetic progression. Iteratively
applying Theorem \ref{lemma:1.4} to
$$A_h+A_{h-1},(A_h+A_{h-1})+A_{h-2},\ldots,\Sum{i=p-2h+2}{h}A_i+A_{p-2h+1},$$
we conclude, in view of (\ref{notallp}) and Corollary
\ref{vosper-cor}, that \be\label{ut} |\Sum{i=p-2h+1}{h}A_i|\geq
3h-p+1,\ee with equality possible only if $\Sum{i=p-2h+1}{h}A_i$ is
an arithmetic progression.

Since $|A_i|\geq 4$ for $i\leq p-2h$, with $|A_i|=4$ possible only
if $A_i$ is not an arithmetic progression, then iteratively applying
Theorems \ref{lemma:1.4} and \ref{vosper-lemma:1.4} to
$$\Sum{i=p-2h+1}{h}A_i+A_{p-2h},\Sum{i=p-2h}{h}A_i+A_{p-2h-1},\ldots,\Sum{i=2}{h}A_i+A_{1},$$
yields, in view of (\ref{notallp}), (\ref{notallp-partii}) and
(\ref{ut}) (note in the last application we may be forced to apply
Theorem \ref{lemma:1.4} instead of Theorem \ref{vosper-lemma:1.4}
even if $|A_h|=4$), that
$$|\Sum{i=1}{h}A_i|\geq 3h-p+1+4(p-2h)-1=3p-5h.$$ Thus
(\ref{notallp}) implies that $h\geq \frac{2p+1}{5}$, whence
$h=\frac{2p+1}{5}$. Hence (\ref{ut}) gives
$|\Sum{i=p-2h+1}{h}A_i|\geq 3$.

Consequently, since $|A_i|\geq 4$ for $i\leq p-2h$, with $|A_i|=4$
possible only if $A_i$ is not an arithmetic progression, then
iteratively applying Theorems \ref{lemma:1.4} and
\ref{vosper-lemma:1.4} to
$$\Sum{i=p-2h+1}{h}A_i+A_{p-2h},\Sum{i=p-2h}{h}A_i+A_{p-2h-1},\ldots,\Sum{i=3}{h}A_i+A_{2},$$
yields, in view of (\ref{notallp-partii}), (\ref{ut}) and
Proposition \ref{propp}, that \be\label{ut2}|\Sum{i=2}{h}A_i|\geq
3h-p+1+4(p-2h-1)=3p-5h-3,\ee with equality possible only if
$\Sum{i=2}{h}A_i$ is an arithmetic progression. Note $h<\frac{p}{2}$
implies that $|A_1| \geq 4$ with equality possible only if $A_1$ is
not an arithmetic progression. Thus, if equality holds in
(\ref{ut2}), then (\ref{notallp}), Theorem \ref{lemma:1.4} and
Corollary \ref{vosper-cor} imply that
\be\label{yet}|\Sum{i=1}{h}A_i|\geq 3p-5h-3+4=3p-5h+1;\ee while on
the other hand, if the inequality in (\ref{ut2}) is strict, then
(\ref{yet}) follows from Theorem \ref{lemma:1.4} and
(\ref{notallp}). Therefore we may assume (\ref{yet}) holds
regardless, whence (\ref{notallp}) implies that $h\geq
\frac{2p+2}{5}$, a contradiction. This completes the proof for
$p>7$.

Suppose $p=7$ and that $\max\{h(S_1), h(S_2)\}\leq 2$. Let
$S_1=U_1U_2$ and $S_2=V_1V_2$ be factorizations of $S_1$ and $S_2$
into square-free subsequences $U_i$ and $V_i$ such that
$|U_2|=|V_2|=3$ and $|U_1|=|V_1|=4$ (as before, such factorizations
are easily seen to exist in view of $\max\{h(S_1), h(S_2)\}\leq 2$
and $|S_1|=|S_2|=p=7$). Let $A_i=U_i\cdot V_i$, and note in view of
Lemma \ref{lemma:1.5} that $|A_i|\geq 4$ for $i=1,2$. Thus applying
Theorem \ref{lemma:1.4} to $A_1+A_2$ implies $|A_1+A_2|=7=p$, so
that the proof is complete in view of $A_1+A_2\subseteq S_1\cdot
S_2$. The case $p=5$ follows by a near identical argument,
concluding the proof.\qfd\end{proof}

\end{document}